\documentclass{amsart}
\usepackage{amsmath,amsthm,amssymb,xypic,array}

\theoremstyle{plain}  

\newtheorem{theorem}{Theorem}[section]
\newtheorem{theoremn}{Theorem}
\newtheorem{proposition}[theorem]{Proposition}
\newtheorem{lemma}[theorem]{Lemma}   
\newtheorem{corollary}[theorem]{Corollary}

\theoremstyle{definition}

\newtheorem{definition}[theorem]{Definition}

\theoremstyle{remark}

\newtheorem{claim}{Claim}

\newtheorem{remark}[theorem]{Remark}

\DeclareMathOperator{\Bsl}{Bsl}

\DeclareMathOperator{\cod}{cod}

\DeclareMathOperator{\Sing}{Sing}

\DeclareMathOperator{\fr}{Frac}
\DeclareMathOperator{\Sec}{Sec}

\newcommand{\QED}{\ifhmode\unskip\nobreak\fi\quad {\rm Q.E.D.}} 

\newcommand\frup[1]{\lceil{#1}\rceil}

\newcommand\kdn[3]{\lceil\frac{{#1\choose #2}}{#3}\rceil}
\newcommand\bin[2]{{#1\choose #2}}

\newcommand{\f}{\varphi}
  
\newcommand{\C}{\mathbb{C}}

\newcommand{\G}{\mathcal{G}}
\newcommand{\Gr}{{G}}
\renewcommand{\H}{\mathcal{H}} 

\newcommand{\I}{\mathcal{I}}

\renewcommand{\L}{\mathcal{L}}

\renewcommand{\O}{\mathcal{O}} 
\renewcommand{\P}{\mathbb{P}}

\newcommand{\T}{\mathbb{T}}

\newcommand{\rat}{\dasharrow}
\renewcommand{\sec}{\mathbb{S}ec}

\newcommand{\exf}{{expected and effective}}
\newcommand{\win}{\sl {\bf wIn}}

\begin{document}

\title{Singularities of linear systems and the Waring Problem}

\author{
Massimiliano Mella}
\address{
Dipartimento di Matematica\\ 
Universit\`a di Ferrara\\
Via Machiavelli 35\\
44100 Ferrara Italia}
\email{mll@unife.it}

\date{April 2004}
\thanks{Partially supported by Progetto Cofin 2002 ``Geometria
  sulle variet\`a algebriche'' MIUR, EAGER}
\maketitle

\section*{Introduction}
Edward Waring stated in 1770 that every integer is a sum of at most 9
positive integral cubes, also a sum of at most 19 biquadrates and so
on. Later on Jacobi and others considered the problem to find all the
decompositions of a given number into the least number of powers,
\cite{Di}.
In this paper I am concerned with  a similar question for general
homogeneous forms. 

Let $f$ be any form of degree $d$ in $n+1$ variables over the
complex numbers. It is easy to prove that $f$ can be decomposed
additively in powers of linear forms. That is
$f=l_1^d+\ldots+l_s^d$ for $l_i$ linear forms. Indeed this is equivalent
to the non degeneracy of the Veronese embedding
$\nu_d(\P^n)\subset\P(S_d)$, where $S_d$ is the vector space of degree
$d$ polynomials in $\C[x_0,\ldots,x_n]$. Moreover the minimum  $s$ for which
$\sec_{s-1}(\nu_d(\P^n))=\P(S_d)$ gives the least number of factors in the
decomposition of a general form $f$. It is now well known that
Alexander--Hirschowitz result, \cite{AH2},
 allows to give a statement, similar to Waring's one, for
addictive decomposition of general homogeneous forms. This has been
first noticed by Iarrobino, \cite{Ia}, using apolarity
and inverse systems.  The way I prefer to look at it, it is through Terracini lemma.

As Jacobi did for the original Waring problem, also the Waring problem
for forms has been  investigated in the
attempt to understand all possible decompositions of a given
form. By means of Hilbert schemes or Grassmannians one can give the structure of an algebraic
variety to the set of decompositions of a given form. In this
way Mukai gave a description of the Fano 3-fold $V_{22}$ as the
variety representing all decompositions of a plane quartic,
\cite{Mu}. Iliev, Ranestad and Schreyer applied similar
arguments to other special cases, \cite{RS}, \cite{IR1} and \cite{IR2}. Classic works of Sylvester,
Reye, Richmond, Hilbert and Palatini investigated the case in which the
decomposition is unique. More recently Iarrobino and Kanev treated the
same problem for special forms, \cite{IK}, see remark \ref{rem:ia}.
 The following list is taken almost verbatim from \cite{RS}.
A general form $f$ of degree $d$ in $n + 1$ variables has a unique
presentation as a sum of 
$s$  powers of linear forms in the following cases:
\begin{itemize}
\item $n = 1$, $d = 2k - 1$ and $s = k$, \cite{Sy}
\item $ n = 3$, $d = 3$ and $s = 5$ Sylvester's Pentahedral Theorem
  \cite{Sy}
\item $n = 2$, $d = 5$ and $s = 7$ \cite{Hi}, \cite{Ri}, \cite{Pa}.
\end{itemize}
In this paper I prove that if $d>n$ there are no further
examples. More precisely I have the following.
\begin{theoremn}\label{th:1}Let $f$ be a general homogeneous form of degree $d$ in
  $n+1$ variables. Assume that $d>n>1$. Then $f$ is
  expressible as sum of $d$th powers of linear forms in a
  unique way if and only if 
  $n=2$ and $d=5$. This unique exception is well known classically,
  \cite{Hi}.
\end{theoremn}
To prove Theorem \ref{th:1} I use the geometric interpretation sketched
before.  In particular I translate the uniqueness
assumption into a statement about singularities of  linear systems with
imposed singularities. Next I study the singularities of this linear
systems. This is an interesting problem in itself. Let me introduce it.
Consider integers $(d,n,l)$. One wants to understand the
singularities of a general hypersurface $F$ of degree $d$ that
is singular at $l$ general points of the projective space $\P^n$. I
expect that when $l$ is as high as possible there is no room, in
general, for other singularities. That is $F$ has only ordinary double
points at the imposed singularities and it is smooth elsewhere. On the
other hand if this is true for the highest $l$ then it is a fortiori
true for lower impositions. Therefore one can expect that, with few
exceptions,  $F$ is irreducible and it has only $l$ ordinary double
points. For $\P^2$ this has been classically studied and it is the
content of a theorem of
Arbarello--Cornalba, \cite{AC}.

The starting point, to extend Arbarello--Cornalba's result to
higher dimension,  are the beautiful papers of
Alexander--Hirschowitz, \cite{AH1} \cite{AH2}.
 I will base not only on their result but
also on their degeneration idea to reach my goal. Indeed to determine the
singularities of $F$ one can consider a suitable degeneration of it. 
This is done, in the full
spirit of A--H, specializing a bunch of
points to an hyperplane. Then an infinitesimal Bertini Theorem,
\cite{CC}, similar to A--C argument for plane curves allows to
determine the singularities on $F$.

The A--H technique consists in choosing
numbers in a clever way. In such a way that one can prove the
statement by a double induction on degree and dimension. Furthermore
when the numbers are not amicable they use the differentiable Horace
Lemma to pop out some points from the hyperplane back to $\P^n$.
 Unfortunately I am not able to apply  Horace Lemma, because I loose control on the
singularities of the specialized linear system. Furthermore
singularities behaves better in a degree induction than in the
dimensional one. For these reasons  I have to develop a slightly different
approach. I play the full induction on the degree when $l$ is
``small'', see Theorem \ref{th:can}. Then
the result I need for Theorem \ref{th:1}, and also a bit more, is obtained in one step using the small $l$ case for degree
$d-1$, see Theorem
\ref{th:fc}. The latter has an application to weakly
defectiveness of Veronese embedding, Corollary \ref{co:ver}.

Here is an outline of the paper. In section \ref{sec:not} I
introduce the main notations, definitions and results that will be used
throughout the paper. In particular a Lemma on the upper semi-continuity
of the dimension of singularities is established. In section
\ref{sec:war}  I explain how to
pass from uniqueness in the Waring problem to a singularity
statement. 
In doing this I show the rationality of
varieties with one apparent $(k+1)$-tuple secant $\P^k$, improving on results of
\cite{CMR}, see also \cite{CR}.

In section \ref{sec:win} I prove
the main numerology to speed up the induction of the main
theorem about singularities of linear systems.   Section
\ref{sec:main}  is the core of the induction and
contains the main result on singularities. Finally I prove Theorem
\ref{th:1}.

I started to work on this problem during a short visit at Universit\`a
Tor Vergata in Roma. It is a pleasure to thank Ciro Ciliberto and
Francesco Russo for discussions, motivations  and  interest in
this project.

\section{Notations and preliminaries}
\label{sec:not}

Unless otherwise stated I work over the field of complex
numbers. First I introduce what is needed to study
linear systems with prescribed singularities.
\begin{definition} Let $p\in\P^n$ be a point. The double point at $p$
  in $\P^n$ is the scheme given by the square of the ideal sheaf of
  $p$. If $P\subset\P^n$ is a collection of points, I denote by $P^2$
  the double points supported on $P$. In particular the linear system 
$|\I_{P^2}(d)|$ is given by hypersurfaces of degree $d$
  singular at $P$.
\end{definition}

Given a collection of points $X=P\cup Q$, with $Q$ supported
on an hyperplane $H$, let $\tilde{X}$ be the residual of $X^2$ with respect
to $H$. 
That is $\tilde{X}=P^2\cup Q$. Then there is the Castelnuovo exact
sequence given by
$$0\to\I_{\tilde{X}}(d-1)\to\I_{X^2}(d)\to\I_{Q^2,H}(d)\to 0$$
This gives the following sequence on cohomology
\begin{equation}
\label{eq:exseq}
0\to H^0(\P^n,\I_{\tilde{X}}(d-1))\to H^0(\P^n,\I_{X^2}(d))\to H^0(H,\I_{Q^2}(d))
\end{equation}

\begin{definition} Consider a collection $P$, of $l$ general points in
  $\P^n$.
Define 
$$\G_{d,n,l}:=|\I_{P^2}(d)|$$
Fix an hyperplane $H$ and a collection $X=P\cup Q$, where $P$ is given by $(l-h)$ general points in
$\P^n$, and $Q=\cup_1^h q_j$ is given by $h$ general points in $H$.
Define
$$\H_{H,d,n,l,h}:=|\I_{X^2}(d)|$$
\label{def:h}

\end{definition}

In this paper I am interested in non empty linear systems of type
$\G_{d,n,l}$, for this I introduce the following definition.

\begin{definition}
 I say that the linear system $\G_{d,n,l}$ is  {\sl expected} if
$$\dim\G_{d,n,l}=\bin{n+d}{n}-(n+1)l-1$$
Moreover if $\G_{d,n,l}$ is expected and $\dim\G_{d,n,l}\geq 0$ I say
that it is 
{\sl expected and effective} 

Note that  if
$\G_{d,n,l}$ is {\exf} then $\G_{d,n,l'}$ is {\exf} for any $l'<l$.
Similarly for linear systems of type $\H$. I say that $\H_{H,d,n,l,h}$ is  {\sl expected and
    effective} 
 if 
$$\dim\H_{H,d,n,l,h}=\bin{n+d}{n}-(n+1)l-1\geq 0$$
Note that if $\H_{H,d,n,l,h}$ is {\exf} then by semi-continuity
$\G_{d,n,l}$ is {\exf}.
In the following I frequently ask {\exf} linear systems of type $\H$ to
satisfy the following further properties. The linear system 
$\H_{H,d,n,l,h}$ is what I need ({\win}) if 
\begin{itemize}
\item $\H_{H,d,n,l,h}$ is {\exf},
\item $|\H_{H,d,n,l,h}\otimes \I_H|\neq \emptyset$
\item   (\ref{eq:exseq}) is a short exact sequence.
\end{itemize}
\end{definition}

The following is a weak form of the main result of \cite{AH2}.

\begin{theorem}[\cite{AH2}] Let $k$ be an
  infinite field and assume that $n\geq 2$ and $d>2$. Then
 $\G_{d,n,l}$ is {\exf} if $l<\kdn{d+n}{n}{n+1}$.

 Furthermore if $\dim\G_{d,n,l}\geq 0$ and $\G_{d,n,l}$  is not {\exf}
 then  $l=\kdn{n+d}{n}{n+1}$, $\dim\G_{d,n,l}=0$, and $(d,n,l)=(3,4,7),(4,2,5),(4,3,9),(4,4,14)$.
\label{th:AH}
\end{theorem}

The aim is to prove the statement on singularities by degeneration. To do this I
need to control how the singularities behaves under 
specializations. This is the content of the next Lemma and Corollary.

\begin{lemma}\label{le:dimsing} Let $\Delta$ be a complex disk around
  the origin. 
Consider the product $V=\P^n\times \Delta$,
 with the natural projections, $\pi_1$ and $\pi_2$. Let $V_t=\P^n\times\{t\}$ and
 $\O_V(d)=\pi_1^*(\O_{\P^n}(d))$.
Fix a configuration $p_1,\ldots,p_l$ of $l$ points on $V_0$ and let
 $\sigma_i:\Delta\to V$ be sections such that $\sigma_i(0)=p_i$ and
 $\{\sigma_i(t)\}_{i=1,\ldots,l}$ are general points of $V_t$ for
$t\neq 0$. Let $P_t=\cup_{i=1}^l\sigma_i(t)$.

 Consider the linear system $\H_t=|\O_{V_t}(d)\otimes\I_{P_t^2}|$. 
Assume that $l<\kdn{d+n}{n}{n+1}$ and $\dim\H_0=\dim\H_t$, for $t\in\Delta$. Let
$\f_i(t):=\dim_{\sigma_i(t)}\Sing\H_t$ and
$\psi_i(t):=\dim_{\sigma_i(t)}\Bsl\H_t$. 
Then for $t\neq 0$,  we have 
$$\f_i(t)\leq \min\{j|\f_j(0)\}$$ 
and 
$$\psi_i(t)\leq \min\{j|\psi_j(0)\}$$
\end{lemma}

\begin{proof} Let $\xi_i$ be the generic point of
  $\sigma_i(\Delta)$ and $P=\cup_{i=1}^l\xi_i$. Let $\H=|\O_{V}(d)\otimes\I_{P^2}|$, then by Theorem \ref{th:AH}, applied to
  the field $\C(t)$, I have
$$\dim\H=\dim\H_t=\bin{n+d}{n}-(n+1)l-1$$
Let $D\in \H$ be a general element then $D_{|V_t}$ is a general
element in $\H_t$. In particular $\Bsl \H_{|V_t}=\Bsl\H_t$ and I can
assume, by Bertini Theorem, that $\Sing\H_t=\Sing\H\cap V_t$, for
$t\neq 0$. The former, together with the
semi-continuity of the dimension of fibers, show that
$\psi_i(t)\leq \psi_i(0)$.
The variety $V_t$ is a Cartier divisor therefore 
$$\dim_{\sigma(0)}(\Sing\H\cap V_0)=\dim_{\sigma(t)}(\Sing\H\cap V_t)$$
This gives
$$\f_i(t)=\dim_{\sigma_i(t)}(\Sing\H\cap V_t)=\dim_{\sigma_i(0)}(\Sing\H\cap
V_0)\leq \dim_{\sigma_i(0)}\Sing\H_0=\f_i(0)$$

To conclude observe that since the points are general a monodromy
argument 
 shows that $\f_i(t)=\f_j(t)$, for $t\neq
0$, and similarly for functions $\psi_i$.
\end{proof}

\begin{corollary}\label{co:sdp} Assume that
  $\H_{H,d,n,l,h}$ is {\win} and there exists $D\in
  \H_{H,d,n,l,h}$ with isolated singularities at some point $p\in P$.
  Then the general element $G\in\G_{d,n,l}$ has
  only ordinary double points.
\end{corollary}

\begin{proof} The linear system $\H_{H,l,d,n,h}$ is {\win}. In
  particular there exists a degeneration like the one described in
  Lemma \ref{le:dimsing} with $\H_0=\H_{H,l,d,n,h}$. The Lemma ensures
  that  $\dim\Sing G=0$. Then by \cite[Theorem 1.4]{CC}
  $G$ has only ordinary double points.
\end{proof}

Let me recall, next, the main definitions and results concerning secant varieties.
Let $\Gr_k=\Gr(k,N)$ be the Grassmannian of $k$-linear spaces in $\P^N$.
Let $X\subset\P^N$ be an irreducible variety
$$\Gamma_k(X)\subset X\times\cdots\times X\times\Gr_k,$$
 the closure of the graph of
$$\alpha:(X\times\cdots\times X)\setminus\Delta\to \Gr_k,$$
taking $(x_0,\ldots,x_{k})$ to the  $[\langle
  x_0,\ldots,x_{k}\rangle]$, for $(k+1)$-tuple of distinct points.
Observe that $\Gamma_k(X)$ is irreducible of dimension $(k+1)n$. 
Let $\pi_2:\Gamma_k(X)\to\Gr_k$ be
the natural projection.
Denote by 
$$S_k(X):=\pi_2(\Gamma_k(X))\subset\Gr_k.$$
Again $S_k(X)$ is irreducible of dimension $(k+1)n$.
Finally let 
$$I_k=\{(x,\Lambda)| x\in \Lambda\}\subset\P^N\times\Gr_k,$$
with natural projections $\pi_i$ onto the factors.
Observe that $\pi_2:I_k\to\Gr_k$ is a $\P^k$-bundle on $\Gr_k$.

\begin{definition} Let $X\subset\P^N$ be an irreducible variety. The {\it abstract $k$-Secant variety} is
$$\Sec_k(X):=\pi_2^{-1}(S_k(X))\subset I_k.$$ While the {\it $k$-Secant variety} is
$$\sec_k(X):=\pi_1(Sec_k(X))\subset\P^N.$$
It is immediate that $\Sec(X)$ is a $((k+1)n+k)$-dimensional variety with a 
$\P^k$-bundle structure on $S_k(X)$. One says that $X$ is
$k$-defective if $$\dim\sec_k(x)<\min\{\dim\Sec_k(X),N\}$$
\end{definition}

\begin{remark} A feature of the above definition is the following
simple observation. Let $\Lambda_1$ and $\Lambda_2$ be two
$(k+1)$-secant $k$-linear space to $X\subset\P^N$. Let $\lambda_1$ and  $\lambda_2$ be
the corresponding projective $k$-spaces in $\Sec(X)$. 
Then we have $\lambda_1\cap \lambda_2=\emptyset$.
\label{re:vuoto} 
\end{remark}
   
Here is the main result I use about secant varieties.
\begin{theorem}[Terracini Lemma \cite{Te}\cite{CC}] \label{th:terracini}
Let $X\subset\P^N$ be an irreducible, projective variety. 
If $p_0,\ldots, p_{k}\in X$ are general points and $z\in\langle
p_0,\ldots, p_{k}\rangle$
 is a general point, then the embedded tangent space at $z$ is
$$\T_z\sec_k(X) = \langle \T_{p_0}X,\ldots, \T_{p_{k}}X\rangle$$
If $X$ is k-defective, then the general hyperplane $H$ containing
$\T_{z}\sec(X)$ is tangent to $X$ along a variety
$\Sigma(p_0,\ldots, p_{k})$ of pure, positive dimension, containing
$p_0,\ldots, p_{k}$.
\end{theorem}

The final statement suggests the following definition

\begin{definition}[\cite{CC}] A variety $X\subset\P^N$ is $k$-weakly
  defective if the general hyperplane, tangent at $(k+1)$ general
  points of $X$, is tangent along a
  positive dimensional subvariety through the tangent points.
\label{def:wd}
\end{definition}

\begin{remark} Terracini Lemma shows that $k$-defective varieties are
  also $k$-weakly defective. The converse is not true in general but
  weak defect is useful to study defective varieties, see
  \cite{CC}.
Note that if $X$ is not weakly defective then the hyperplane
is tangent exactly at the $(k+1)$ points and has ordinary double
points there, \cite[Theorem 1.4]{CC}.
\end{remark}

\section{From uniqueness to singularities}
\label{sec:war}

 Let $\nu_d(\P^n)\subset\P^N$ be the $d$-uple Veronese embedding of
$\P^n$. I already observed, see
\cite{Ci} for a deeper  account,
that the problem of finding the minimal number of linear forms needed to
express a general form of degree $d$ is equivalent to determine the
dimension of $\sec_k(\nu_d(\P^n))$. 
A general form $f$ of degree $d$ is just a general point $p\in\P^N$. The number of
representations  of $f$ as a sum of $(k+1)$ powers of
linear forms correspond to the number of $(k+1)$-secant linear spaces to
$\nu_d(\P^n)$ passing through $p$. If there is just one
I will say that there exists the {\sl canonical form}.
Note that $N=\bin{d+n}{n}-1$ and $\dim\Sec_{k}(\nu_d(\P^n))=(n+1)(k+1)-1$
therefore there are finitely many representations only if
$k+1=\frac{\bin{d+n}{n}}{n+1}$ is an integer.

 Assume that $k+1=\frac{\bin{d+n}{n}}{n+1}$ is an integer. 
Consider the natural
map 
$$\pi_1:\Sec_{k}(\nu_d(\P^n))\to\P^n$$
 then the canonical
form exists only if $\pi_1$ is dominant and  birational.
The following is a generalization of \cite[Theorem
  4.1]{CMR}, see also \cite[Theorem 2.6]{CR}.

\begin{theorem} Let $X\subset\P^N$ be an irreducible variety of dimension $n$. 
Assume that the natural map $\sigma:\Sec_{k}(X)\to\P^N$ is dominant and
birational. Let $z\in\sec_{k-1}(X)$ be a general point.  
Consider $\f:\P^N\rat\P^n$ the projection from the embedded tangent
space $\T_z\sec_{k-1}(X)$. Then
$\f_{|X}:X\rat \P^n$ is dominant and birational.
\label{th:proj}
\end{theorem}

\begin{proof} Choose a general point $z$ on a general $k$-secant $\P^{k}=\langle
  p_0,\ldots,p_{k-1}\rangle$.
Let $f:Y\to \P^N$ be the blow up of $\sec_{k-1}(X)$ with exceptional
divisor $E$, and fiber  $F_z=f^{-1}(z)$. 
Let $y\in F_z$ be a general point. This point uniquely determines a
linear space $\Pi$  of dimension $k(n+1)$ that
contains $\T_z\sec_{k-1}(X)$. 
Then the projection $\f_{|X}:X\rat \P^n$ is birational if and only if
$(\Pi\setminus \T_z\sec_{k-1}(X))\cap X$ consists
of just one point.

Assume that there exist two points $x_1$ and $x_2$ in $(\Pi\setminus
\T_z\sec_{k-1}(X))\cap X$.
By Terracini Lemma, Theorem \ref{th:terracini},
$$\T_z\sec_{k-1}(X)= \langle \T_{p_0}X,\ldots, \T_{p_{k-1}}X\rangle$$
Consider the linear spaces $\Lambda_1=\langle x_1, p_0,\ldots,p_{k-1}\rangle$ and 
$\Lambda_2=\langle x_2,p_0,\ldots,p_{k-1} \rangle$. The Trisecant
Lemma, see for instance \cite[Proposition 2.6]{CC}, yields
$\Lambda_1\neq\Lambda_2$. Let
$\Lambda_1^Y$, $\Lambda_2^Y$ and $\Pi^Y$ be the strict transforms on $Y$. 
Since $z\in\langle p_0,\ldots,p_{k-1}\rangle$ and $y=\Pi^Y\cap F_z$ then both $\Lambda_1^Y$ and
$\Lambda_2^Y$ contain the point $y\in F_z$. In particular I have
$$\Lambda_1^Y\cap\Lambda_2^Y\neq\emptyset$$ 

Let $\pi_1:\Sec_k(X)\to\P^N$ be the morphism from the abstract secant variety,
and  $\mu:U\to Y$ the induced morphism. That is
$U=\Sec_k(X)\times_{\P^N} Y$. Then 
there exists a commutative diagram
$$\diagram
U\dto_{p}\rto^{\mu}&Y\dto^{f}\\
\Sec_k(X)\rto^{\pi_1}&\P^N\enddiagram$$

Let $\lambda_i$ and $\Lambda^U_i$ be the strict transform of $\Lambda_i$ in $\Sec_k(X)$ and $U$ 
respectively. By Remark \ref{re:vuoto}
 $\lambda_1\cap \lambda_2=\emptyset$, so that
$$\Lambda^U_1\cap \Lambda^U_2=\emptyset.$$

This proves that $\mu^{-1}$ is not defined  on $y$, a general point of
the divisor
$E$. That is $\mu$, and henceforth $\pi_1$, are not birational.
\end{proof}

\begin{remark} A variety $X^n\subset\P^N$ is said to have one apparent
  $(h+1)$-tuple $\P^h$ if there exists and is unique a
  $(h+1)$-secant $\P^h$ through the general point of $\P^N$, for wider
  discussions and properties of these varieties I recommend
  \cite{CR}. In this dictionary Theorem \ref{th:proj} proves that a variety with  one apparent
  $(h+1)$-tuple $\P^h$ is rational. This extends the result of
  \cite{CMR} for varieties with one apparent double point, see also \cite{CR}.
\end{remark}
 
Let  $k+1=\frac{\bin{d+n}{n}}{n+1}$ be an integer and $\f:\P^N\rat \P^n$ the
projection from a  general tangent space to
$\sec_{k-1}(\nu_d(\P^n))$. Theorem \ref{th:proj}, says that the
canonical form exists, only if $\f_{|\nu_d(\P^n)}$ is birational. Terracini Lemma shows that
$\f_{|\nu_d(\P^n)}$ is the map induced by $\G_{d,n,k}$. Therefore to have the
canonical form the linear system $\G_{d,n,k}$ has to be {\exf} and 
the map associated to $\G_{d,n,k}$ has
to be birational onto $\P^n$.
 The final brick of the bridge, is the original statement of
 Noether--Fano inequalities, reinterpreted in modern terminology.

\begin{theorem}[Noether--Fano inequality \cite{Co}] Let
 $\L$ be a linear system, without fixed components, of forms of degree
$d$  on $\P^n$.  
Assume that $\dim\L=n$ and the rational map associated to $\L$,
 is birational. Then $(\P^n,(n+1)/d\L)$ has not canonical
  singularities.
\label{th:NF}
\end{theorem}

I have proved the following necessary condition
\begin{proposition} Let $f$ be a general homogeneous form of degree $d$ in
  $n+1$ variables, with $n\geq 2$.  Then $f$ is
  expressible as sum of $(k+1)$ powers of linear forms in a
  unique way  only if the pair $(\P^n,(n+1)/d\G_{d,n,k})$ has not
  canonical singularities, where $k=\frac{\bin{d+n}{n}}{n+1}-1$ is an integer.
\label{pro:bridge}
\end{proposition}

To test this condition I have to understand the singularities of the
linear system $\G_{d,n,k}$.  I will do it by a degeneration argument. First I have to
establish some good degenerations. This is the content of the following
section.

\section{Numerical conditions to be {\win}}
\label{sec:win}

In this section I establish sufficient conditions for a linear system of
type $\H$ to be {\win}.   This is just
a,  maybe pedant, way to break the
main induction argument in small pieces.

\begin{lemma} Let $H\subset\P^n$ be an hyperplane and  
assume that $\G_{d,n,l}$, $\G_{d-1,n,l-h}$, are {\exf} and
 $\G_{d,n-1,h}$ is expected.
 Assume further that
\begin{equation}
\label{eq:bound}
\bin{n+d-1}{n}-(n+1)(l-h)-h\geq\max\{1,\dim\G_{d-2,n,l-h}\}
\end{equation}
 Then $\H_{H,d,n,l,h}$ is {\win}.
\label{le:dimbase}
\end{lemma}
\begin{proof} Let $X=P\cup Q$ be as in Definition \ref{def:h}, where
  $Q=\cup_1^h q_j\subset H$.
By hypothesis $\G_{d,n,l}$ is {\exf} and by semi-continuity I have  to
prove that 
$$\dim\H_{H,d,n,l,h}+1\leq \bin{n+d}{n}-(n+1)l$$
Let 
$$h(d-1)=\dim H^0(\P^n,\O_{\P^n}(d-1)\otimes \I_{\tilde{X}})$$
and
$$h(n-1)=\dim H^0(\P^{n-1},\O_{\P^{n-1}}(d)\otimes \I_{Q^2})$$
By sequence (\ref{eq:exseq}) I have
\begin{equation}
\label{eq:findim}
\dim\H_{H,d,n,l,h}+1\leq h(d-1)+h(n-1)
\end{equation}
By hypothesis  $\G_{d,n-1,h}$ is
expected therefore 
$$h(n-1)=\bin{n-1+d}{n-1}-nh$$
The claim is therefore equivalent to prove that 
\begin{equation}
\label{eq:tbp}
h(d-1)=\bin{n+d-1}{n}-(n+1)(l-h)-h>0
\end{equation}
Indeed this yields 
$$\dim\H_{H,d,n,l,h}+1\leq \bin{n+d-1}{n}+\bin{n-1+d}{n-1}-(n+1)l=\bin{n+d}{n}-(n+1)l$$
and $\dim|\H_{H,d,n,l,h}\otimes\I_H|=h(d-1)-1\geq 0$. As a consequence
I get the
exactness of sequence (\ref{eq:exseq}).

By hypothesis $\G_{d-1,n,l-h}$ is {\exf}, that is 
$$\dim\G_{d-1,n,l-h}=\bin{n+d-1}{n}-(n+1)(l-h)\geq 0$$
Hence to prove equation
(\ref{eq:tbp}) it is enough to show that the $h$ points on $H$
impose independent conditions to $\G_{d-1,n,l-h}$. This is equivalent
to ask that
$$H\not\subset\Bsl|\I_{P^2}(d-1)\otimes 
(\otimes_{j=1}^{h-1}\I_{q_j})|$$
(the $h-1$ is not a misprint).
That is to say
$$\bin{n+d-1}{n}-(n+1)(l-h)-(h-1)>\dim\G_{d-2,n,l-h}$$
I can therefore conclude by equation (\ref{eq:bound}) that 
$$h(d-1)=\bin{n+d-1}{n}-(n+1)(l-h)-h>0$$
\end{proof}

Next I translate the {\exf} conditions of Lemma \ref{le:dimbase} into
numerical conditions.

\begin{lemma} Fix  integers $n\geq 3$, $d\geq 3$, $l>h$,  and an hyperplane
  $H\subset\P^n$, together with the following conditions:
\begin{itemize}
\item[(L)] $l<\kdn{n+d}{n}{n+1}$
\item[(H)] either $h<\kdn{n-1+d}{n-1}{n}$ or
  $h=\frac{\bin{n-1+d}{n-1}}{n}$ and $\G_{d,n-1,h}$ is expected
\item[(LH)] $l-h<\kdn{n+d-1}{n}{n+1}$
\item[(C)] $\bin{n+d-1}{n}-(n+1)(l-h)-h> 0$
\item[(D4)] $l-h> n$
\item[(D3)] $l<\frac{\bin{n+3}{n}}{n+1}-\frac{n+2}3+1$ and
$h=l-1$.
\end{itemize}

Then the linear system $\H_{H,d,n,l,h}$ is {\win} if one of the
following set of conditions are satisfied:
\begin{itemize}
\item $d\geq 5$ and (L), (H), (LH), (C)
\item $d=4$ and (L), (H), (LH), (C), (D4)
\item $d=3$ and (D3).
\label{le:dim}
\end{itemize}
\end{lemma}
\begin{proof} I will check that the hypothesis of Lemma
  \ref{le:dimbase} are satisfied for every set of conditions.
Hypothesis (L) and  (H) together with Theorem \ref{th:AH} ensure that
$\G_{d,n,l}$ is {\exf} and
$\G_{d,n-1,h}$ is expected.

First note that for $d=3$ condition (D3) forces condition (L), (H), and (C).
Indeed I have
$$h<\frac{\bin{n+3}{n}}{n+1}-\frac{n+2}3=\frac{(n+3)(n+2)-2(n+2)}{6}=\frac{\bin{n+2}{n-1}}{n}$$
and
$$\bin{n+2}{n}-(n+1)-\frac{\bin{n+3}{n}}{n+1}+\frac{n+2}3\geq
\frac{2(n+2)(n+1)}6-(n+1)> 0$$
Moreover  $\G_{2,n,1}$ is {\exf} and $\G_{1,n,1}$ is empty. Therefore
I can apply Lemma \ref{le:dimbase} to conclude.

If $d>3$ condition (LH), together with Theorem \ref{th:AH} shows that $\G_{d-1,n,l-h}$ is {\exf}.

If $\G_{d-2,n,l-h}$ is empty then, 
by condition (C), I conclude by Lemma \ref{le:dimbase}. This is always
the case for $d=4$.

Assume that  $d\geq5$ and $\dim\G_{d-2,n,l-h}$ is not empty. If
$\G_{d-2,n,h}$ is not {\exf} then by Theorem \ref{th:AH} 
 $\dim\G_{d-2,n,l-h}=0$ and by condition (C) I conclude with Lemma
\ref{le:dimbase}. 
If $\G_{d-2,n,h}$ is {\exf} then note that
$$h\leq \kdn{n-1+d}{n-1}{n}-1< \bin{n+d-2}{n-1}$$
therefore condition  (H) yields
$$\dim\G_{d-2,n,l-h}= \bin{n+d-2}{n}-(n+1)(l-h)-1<
\bin{n+d-1}{n}-(n+1)(l-h)-h-1$$
and I can apply Lemma \ref{le:dimbase} to conclude.
\end{proof}

The following Lemma allows to play induction for
``small'' $l$.
Before starting the proof it is convenient to introduce a definition

\begin{definition} For integers $a$ and $b$ let
$\frup{a,b}:=\kdn{a+b}{a}{a+1}-\frac{\bin{a+b}{a}}{a+1}$.
\end{definition}

\begin{lemma} Fix an hyperplane $H\subset\P^n$. Let 
$$l_d:=l=\kdn{n+d+1}{n}{n+1}-\kdn{n+d}{n-1}{n}$$
and 
$$h_d:=h=\kdn{n+d+1}{n}{n+1}-\kdn{n+d}{n-1}{n}-\kdn{n+d}{n}{n+1}+\kdn{n-1+d}{n-1}{n}$$ 
Assume that $d\geq 4$ and $n\geq 3$ 
then $\H_{H,d,n,l,h}$ is {\win}.
\label{le:lh}
\end{lemma}
\begin{proof} Note that $(l_d-h_d)=l_{d-1}$.
First I simplify the expressions of $h$ and $l$,
\begin{eqnarray}
h&=&\kdn{n+d+1}{n}{n+1}-\kdn{n+d}{n-1}{n}-\kdn{n+d}{n}{n+1}+\kdn{n-1+d}{n-1}{n}\nonumber\\
&=&\frac{n\bin{n+d}{n-1}-(n+1)\bin{n+d}{n-1}+(n+1)\bin{n-1+d}{n-1}}{n(n+1)}
+\fr(h)\nonumber\\
&=&\frac{\bin{n-1+d}{n-1}}n-\frac{\bin{n+d}{n-1}}{n(n+1)}+ \fr(h)
\label{eq:h}
\end{eqnarray}
where I collected all the fractional parts in the term 
\begin{equation}
\label{eq:frh}
\fr(h)=\frup{n,d+1}-\frup{n-1,d+1}-\frup{n,d}+\frup{n-1,d}
\end{equation}
in particular $-2<\fr(h)<2$;

\begin{eqnarray}
l&=&\kdn{n+d+1}{n}{n+1}-\kdn{n+d}{n-1}{n}=
\frac{n\bin{n+d+1}{n}-(n+1)\bin{n+d}{n-1}}{n(n+1)}+\fr(l)\nonumber\\
&=&\frac{\bin{n+d}{n}}{n+1}-\frac{\bin{n+d}{n-1}}{n(n+1)}+\fr(l)
\label{eq:l}
\end{eqnarray}
where the fractional part is
\begin{equation}
\label{frl}
\fr(l)=\frup{n,d+1}-\frup{n-1,d+1}
\end{equation}
this time $-1<\fr(l)<1$. Incidentally this shows that both $l$ and
$l-h$ are positive integers.

The plan is to prove the claim by checking the
conditions of Lemma \ref{le:dim}.

Condition (H) in Lemma \ref{le:dim} is given by equation (\ref{eq:h}) and 
 the following inequality
\begin{equation}
\label{ineq:h}
\bin{n+d}{n-1}\geq 2n(n+1)
\end{equation}
which is true for any $d\geq 4$ and $n\geq 4$. In the missing case,
%
%
namely $(n,d)=(3,4)$, I have $h=3<5=\kdn{6}{2}{3}$.

Condition (L) is obtained by equation (\ref{eq:l}) and the following inequality
$$\bin{n+d}{n-1}\geq n(n+1)$$
which is verified for any $d\geq 3$ and $n\geq 3$.
%
%

As observed condition (LH) is just condition (L) for $d-1$.

To check condition (C) note that by equations (\ref{eq:h}) and (\ref{eq:l})
\begin{eqnarray*}
\bin{n+d-1}{n}-(n+1)(l-h)-h&=&\bin{n+d-1}{n}-(n+1)l+nh\\
&=&\frac{\bin{n+d}{n-1}}{n}-\frac{\bin{n+d}{n-1}}{n+1}-(n+1)\fr(l)+n\fr(h)\\
&\geq&\frac{\bin{n+d}{n-1}}{n(n+1)}-(n\frup{n,d}+\frup{n,d+1})
\end{eqnarray*}
The requested bound is therefore implied by
$$\bin{n+d}{n-1}\geq n(n+1)^2$$
This is verified for $d\geq 4$ and $n\geq 7$ or $d\geq 5$ and $n\geq
5$ or $d\geq 6$ and $n\geq 4$ or $d\geq 8$ and $n\geq 3$. 
%
%
The remaining cases can be checked by a direct computation.
Let $\delta=\bin{n+d-1}{n}-(n+1)(l-h)-h$ then
$${\setlength{\extrarowheight}{4pt}
\begin{tabular}{|c|cccccccc|}
\hline
$(d,n)$&(4,6)&(4,5)&(4,4)&(4,3)&(5,4)&(5,3)&(6,3)&(7,3)\\
\hline
$l-h$&9&7&5&4&12&7&11&18\\
\hline
$h$&15&9&7&3&9&4&7&9\\
\hline\hline
$\delta$&6&5&3&1&1&3&5&3\\
\hline
\end{tabular}}
$$

If $d=4$ I have to prove (D4),
by definition
$$l-h=\kdn{n+4}{n}{n+1}-\kdn{n+3}{n-1}{n}>\frac{\bin{n+3}{n}}{n+1}-\frac{\bin{n+3}{n-1}}{n(n+1)}-1$$
Expanding the binomial I get
$$l-h>\frac{(n+3)(n+2)}{8}-1$$
The requested bound is therefore implied by
$$(n+3)(n+2)\geq 8(n+1)$$
This is verified for $n\geq 4$. For $n=3$ I have $l-h=9-5>3=n$.
\end{proof}
\section{Singularities of linear systems}
\label{sec:main}

 This is the unique result I have about singularities of cubic forms.
\begin{theorem} Let $G\in\G_{d,n,l}$ be a general element. Assume that
either
\begin{itemize}
\item  $d\geq 4$, $n\geq 3$, and
  $l=\kdn{n+d+1}{n}{n+1}-\kdn{n+d}{n-1}{n}$, or
\item $d=3$, $n\geq 3$ and $l<\frac{\bin{n+3}{n}}{n+1}-\frac{n+2}3+1$.
\end{itemize}
 Then $\G_{d,n,l}$ is {\exf}  and $G$ has
 only ordinary double points.
\label{th:can}
\end{theorem}

\begin{remark} Like
  in the interpolation problem cubics are kind of difficult because
  the degeneration technique does not  work for the highest possible values of $l$.
Nevertheless this is what is needed to conclude  for degree
  at least 4.
\end{remark}

\begin{proof}
I prove the claim by induction on $d$. 

Assume that $d\geq 4$.
Let $l'=\kdn{n+d}{n}{n+1}-\kdn{n+d-1}{n-1}{n}$ and 
$h=l-l'$. Fix a general hyperplane $H\subset\P^n$, then 
by Lemma \ref{le:lh} $\H_{H,d,n,l,h}$ is {\win}. This proves that
$\G_{d,n,l}$ is {\exf}. Since $|\H_{H,d,n,l,h}\otimes\I_H|\neq \emptyset$ I can look
 the linear system $\H_{H,d,n,l,h}$ from a different viewpoint.
Choose $l-h$ general points in $\P^n$, $p_1,\ldots,p_{l-h}$.
Let 
$$D\in \G_{d-1,n,l-h}$$ 
be a general element.
Choose a general hyperplane $H\subset \P^n$
and $h$ general points on $C=D\cap H$, $q_1,\ldots,q_h$.

By induction hypothesis $D$ has only ordinary double points.
Therefore the general element of
$\H_{H,d,n,l,h}$ has only isolated singularities at the $p_i$s and by
Corollary \ref{co:sdp} the general element of $\G_{d,n,l}$ has only ordinary
double points.

To conclude I need to prove the first step of induction.
 If $d=3$ one gets
$$l=\kdn{n+4}{n}{n+1}-\kdn{n+3}{n-1}{n}=\frac{(n+3)(n+2)}8+\frup{n+4,n}-\frup{n+3,n-1}$$
and 
$$l-\frac{\bin{n+3}{n}}{n+1}+\frac{n+2}3-1<\frac{(n+2)(5-n)}{24}$$
This proves condition (D3) for $n\geq 5$. Therefore $\H_{H,3,l,l-1}$
is {\win} and there is at least one divisor $D\in\H_{H,3,l,l-1}$ with
$D=Q+H$ with $Q$ a quadric cone. This together with Corollary
\ref{co:sdp} proves the theorem for $n\geq 5$. 
For $n=3,4$ note that $l$ is, respectively, at most 4 and 5. Then I do not need any
specialization. The linear system $\G_{d,n,l}$ always contain a
reducible divisor $Q+H$ where $Q$ is a quadric cone.
\end{proof}

I am now ready to prove the main result about forms
with assigned singularities.

\begin{theorem} Assume that $n\geq 3$ and $d\geq 3$ and either
\begin{itemize}
\item $l_0=\kdn{n+d+1}{n}{n+1}-1$ and $n\frup{n-1,d+1}-(n+1)\frup{n,d+1}+1>0$,
  or
\item  $l_1=\kdn{n+d+1}{n}{n+1}-2$ or
\item $l_2=\kdn{n+d+1}{n}{n+1}-1$, $\frup{n-1,d+1}=0$ and either
  $d\geq 4$ or $d=3$ and $n\geq 6$.
\end{itemize}
Then the general element in $\G_{d+1,n,l_i}$ has only ordinary double points.
\label{th:fc}
\end{theorem}

\begin{remark} 
The $n=2$ case is well known, \cite{AC}. I will recall it in Lemma
  \ref{le:p2}. Note that for $n=2$, $d=6$, and $l=9$ the unique 6-ic
  with 9 general double points is a double cubic.
For $n=3$ there is another exception. Consider $\G_{4,3,8}$ a dimension count shows
  that $\G_{4,3,8}=\{Q_1^2,Q_2^2,Q_1Q_2\}$, where $Q_i$ are the
  quadrics passing simply through the 8 points. In particular the
  general element of $\G_{4,3,8}$ has a curve of singularities. I believe that there
  are very few {\exf} linear systems $\G_{d,n,l}$ with non isolated
  singularities. 
Maybe these are the only ones.
 It would be interesting to classify them all, see Corollary
  \ref{co:num} for a partial result in this direction. Note that all
  linear systems that are not expected have non isolated singularities.
\end{remark}
\begin{corollary} The Veronese embedding $\nu_d(\P^n)\subset\P^N$ is
  not $k$-weakly defective for $d\geq 4$ and $n\geq 3$ if
  $\cod\sec_k(\nu_d(\P^n))\geq n+1$.
\label{co:ver}
\end{corollary}

\begin{remark}\label{rem:ia} Luca Chiantini and Ciro Ciliberto pointed me
  out the following consequence of non  weakly defectiveness of the
  Veronese embedding.
Assume that $\nu_d(\P^n)$ is not $k$-weakly defective. Let
  $x\in\nu_d(\P^n)$ be a general point. Then, by Terracini Lemma, there exists a unique
  $\P^k$ which is $(k+1)$-secant through $x$. That is the natural map
  $\pi_k:Sec_k(\nu_d(\P^n))\to \P^N$ is birational onto the
  image. This improves results of \cite[Theorems 2.6, 7.18]{IK} where this statement was
  proved for $k\sim \bin{n+d/2}{n}$.
\end{remark}

\begin{proof}[Proof of the Theorem]
Let
 $h_0=h_1=\kdn{n+d}{n-1}{n}-1$ and
 $h_2=\frac{\bin{n+d}{n-1}}{n}$ be integers. Note that if either
 $d\geq 4$ or $d=3$ and $n\geq 6$ then 
 by Theorem \ref{th:AH} the linear system
 $\G_{d+1,n-1,h_2}$ is expected. Theorem \ref{th:can} shows that  the general element
$D\in\G_{d,n,l_i-h_i}$ has only ordinary double points. Consider the linear systems
$\H_{H,d+1,n,l_i,h_i}$. By hypothesis I have
$$
\bin{n+d}{n}-(n+1)(l_0-h_0)-h_0=n\frup{n-1,d+1}-(n+1)\frup{n,d+1}+1>0
$$
and, for $i=1,2$
$$
\bin{n+d}{n}-(n+1)(l_i-h_i)-h_i=(n+1)+n\frup{n-1,d+1}-(n+1)\frup{n,d+1}+1>0
$$

This proves, via Lemma \ref{le:dim}, that $\H_{H,d+1,n,l_i,h_i}$ is
{\win}, for $d\geq 4$. 
For $d=3$ note that
\begin{eqnarray*}
l_i-h_i&>& \frac{n\bin{n+4}{n}-(n+1)\bin{n+3}{n-1}}{n(n+1)}-2=
\frac{n\bin{n+3}{n}-\bin{n+3}{n-1}}{n(n+1)}-2\\
&=&\frac{(n+3)(n+2)}{8}-2
\end{eqnarray*}
for $i=0,1,2$.
Then $l_i-h_i>n$ for $n\geq 5$ and, via Lemma \ref{le:dim},  $\H_{H,d+1,n,l_i,h_i}$ is
{\win} in these cases.
The leftover has to be checked directly.
If $(d,n)=(3,4)$ I have $l_0-1=l_1=12$ and $h_1=8$. 
If $(d,n)=(3,3)$ I have $l_1=7$, $h_1=3$. Therefore I
can apply Lemma \ref{le:dim} for $\H_{H,4,3,l_i,h_i}$, with
$i=0,1$. and $\H_{H,4,4,l_0,h_0}$. 
For $\H_{H,4,4,l_1,h_1}$ note that 
$$\bin{n+3}{n}-(n+1)(l_1-h_1)-h_1=7>\dim\G_{2,4,4}=1$$
This time I use Lemma \ref{le:dimbase} to conclude that $\H_{H,4,4,l_1,h_1}$ is
{\win}.

I can, therefore, assume that $D+H\in\H_{H,d+1,n,l_i,h_i}$ for $D\in
\G_{d,n,l_i-h_i}$ with only ordinary double points. This
together with  Corollary \ref{co:sdp} conclude the proof.
\end{proof}

I can use the above theorem to determine the singularities of
$\G_{d,n,l}$ for $d\geq 5$ prime and $d=4$.

\begin{corollary} Assume that $\G_{d,n,l}$ is effective. Then the
 general element $G\in\G_{d,n,l}$ has ordinary double points
 if either $d=4$, and $(n,l)\neq (3,8),(3,9)$ or $d\geq 5$ is
  a prime. I already
  described the exceptions.
\label{co:num}
\end{corollary}
\begin{proof}
 Note that
$$\frac{\bin{n+d}{n}}{n+1}=\frac{(n+d)(n+d-1)\cdot\ldots\cdot(n+2)}{d(d-1)!}$$

Note that both
$$\frac{(n+d)(n+d-1)\cdot\ldots\cdot(n+2)}{(d-1)!}$$
and
$$\frac{(n-1+d)(n-1+d-1)\cdot\ldots\cdot(n+1)}{(d-1)!}$$
are integers. Furthermore for some $1\leq a\leq d$ I have $n+a\equiv
0(d)$. Assume that $d$ is a prime. Then either $\frup{n,d}$ or $\frup{n-1,d}$
vanish and I can apply Theorem \ref{th:fc} for any triple $(d,n,l)$
with $\G_{d,n,l}$ effective.

Assume that $d=4$, this time
$$\frac{\bin{n+4}{n}}{n+1}=\frac{(n+4)(n+3)(n+2)}{4\cdot3\cdot2}$$
If $n+2\equiv 0(2)$ I have $(n+4)(n+3)(n+2)\equiv
0(24)$. While if $n+2\equiv 1(2)$ I have $(n+3)(n+2)(n+1)\equiv 0(24)$.
This yields that either $\frup{n,4}=0$ or
$\frup{n-1,4}=0$. This proves the claim, by Theorem \ref{th:fc}, for $n\geq 6$.
If $n=2,4,5$ it is easy to check that $\frup{n,4}=0$. This
concludes the proof.
\end{proof}

\begin{remark}
 Note that for $(d,n,l)=(6,9,500)$ the conditions of Theorem \ref{th:fc}
are not satisfied. This is the first occurrence that I cannot treat in
degree 6. I have not statement, similar to Corollary \ref{co:num}, 
for fixed $n$ and any $d$.
\end{remark}

\section{Proof of Theorem \ref{th:1}}

I am now in the condition to prove Theorem \ref{th:1}.
The following lemma allows to treat the two dimensional  case.

\begin{lemma}[\cite{AC}] Let $l=\kdn{d+2}{2}{3}-1$,
then $(\P^2,1/2\G_{d,n,l})$ has canonical singularities.
\label{le:p2}
\end{lemma}
\begin{proof}Let $G\in\G_{d,n,l}$ be a general element. Then by
  \cite[Theorem 1.4]{CC}
  either $G$ has only the imposed ordinary double points, or it has a fixed
  component of singularities passing through $P$. The latter is easily
  seen, by a dimension count, to be
  possible only for $d=6$ and $l= 9$. In this special case the unique
  6-ic with assigned 9 double points is a double cubic.
Thus $G$ has only double points and $(\P^2,1/2\G_{d,n,l})$ has
canonical singularities.
\end{proof}

\begin{proof}[Proof of Theorem \ref{th:1}] I already observed, 
 proposition \ref{pro:bridge}, that for the existence of the canonical
  form 
$$k+1=\frac{\bin{d+n}{n}}{n+1}$$
 has to be an integer, and $(\P^n,(n+1)/d\G_{d,n,k})$ has to be not
 canonical.

If $n=2$ this means that either $d=4,5$ or $d\geq 7$. If $d=4$ I have
$k=4$, and $d^2=16=4k$. In particular
the map associated to $\G_{4,2,4}$ is composed with a pencil. Note
that $\dim\G_{4,2,5}=0$ and $\G_{4,2,5}$ is not expected. If $d=5$
I have $k=6$ and $d^2=25=4k+1$. It is not difficult to see that
the scheme base locus of $\G_{5,2,6}$ is just $P^2$. In this case we
already known that there exists a canonical form, \cite{Hi}.
Assume that $d\geq 7$ then \cite[Theorem 3.2]{AC}  shows that
$\G_{d,2,k}$ does not have fixed components and by Lemma \ref{le:p2} $(\P^2,3/d\G_{d,2,k})$
  has canonical singularities. Therefore uniqueness is impossible.

If $n\geq 3$ then by Theorem \ref{th:fc} the general element
$G\in\G_{d,n,k}$  has only ordinary double points. In particular the
divisor $G$
is irreducible and  after the blow up of the singular points it is
smooth. This proves that the log pair
 $(\P^n,(n+1)/d\G_{d,n,k})$ has canonical singularities and concludes the proof.
\end{proof}

\begin{remark} To extend Theorem \ref{th:1} to lower degree one has
  to study the base locus of $\G_{d,n,l}$ and not only its
  singularities. The semi-continuity works for the dimension of the
  Base Locus. Furthermore, with similar argument it is possible to prove that functions $\psi_i$
  of Lemma \ref{le:dimsing} are zero at least for ``small'' $l$. What
  is completely  missing is a way to determine the scheme base locus
  out of these informations.
\end{remark}

Theorem \ref{th:1} together with the known exceptions allow to
answer the uniqueness question for forms of at most 4 variables.
\begin{corollary} Let $f$ be an homogeneous form of degree $d$ in $n+1$
  variables.
If $n\leq 3$ then canonical form exists if and only if $(d,n)=(2k+1,1),(5,2),(3,3)$.
\end{corollary}

\end{document}